\documentclass[12pt]{article}
\usepackage{euscript}
\usepackage{amsfonts}
\usepackage{amsthm}
\usepackage{amsmath}
\usepackage{amsfonts}
\usepackage{latexsym}
\usepackage{amssymb}
\usepackage[ansinew]{inputenc}

\begin{document}

\newtheorem{lemma}{Lemma}[section]
\newtheorem{theo}[lemma]{Theorem}
\newtheorem{coro}[lemma]{Corollary}
\newtheorem{rema}[lemma]{Remark}
\newtheorem{pro}[lemma]{Proposition}
\newtheorem{fed}[lemma]{Definition}
\def\bdem{\begin{proof}}
\def\edem{\end{proof}}
\def\bequ{\begin{equation}}
\def\eequ{\end{equation}}
\renewcommand{\thesection}{\arabic{section}}
\renewcommand{\theequation}{\thesection.\arabic{equation}}
\newcommand{\equnew}{\setcounter{equation}{0}}

\newcommand{\laba}{\label}

%%%%%%%%%%%%%%%%%%%%%%%%%%%%%%%%%%%%%%%%%%%%%%%%%%%%%%%%%%%%%%%%%%%%%%%%

\newcommand{\rr}{\mbox{$    \rightarrow   $}}
\newcommand{\disc}{ {\Bbb D} }
\newcommand{\noi}{\noindent}
\newcommand{\ov}{ \overline }
\newcommand{\om}{\omega}
\newcommand{\dist}{\mbox{dist}\,}
\newcommand{\ovp}{\ov{\partial}}
\newcommand{\supp}{ \mbox{supp}\,  }
\newcommand{\Om}{\Omega}
\newcommand{\la}{\langle}
\newcommand{\ra}{\rangle}
\newcommand{\berg}{A^2}
\newcommand{\deltab}{\tilde{\Delta}}
\newcommand{\linf}{L^\infty(\disc)}
\newcommand{\radop}{T_{rad}}
\newcommand{\li}{\ell^{\infty}}
\newcommand{\peso}[1]{ \quad \text{ \rm  #1 } \quad }
\renewcommand{\qed}{\hfill $\blacksquare$}
\newcommand{\eps}{\varepsilon}
\newcommand{\clau}[1]{\overline{#1}^{\ell^\infty}}
\newcommand{\sgn}{\mbox{sgn}\,}
\newcommand{\piso}[1]{\ \ \mbox{#1}\ \ }

%%%%%%%%%%%%%%%                   FONTS
\newcommand{\toep}{\mathfrak{T}}
\newcommand{\cald}{\mathfrak{D}}
\renewcommand{\L}{\mathfrak{L}}
\def\N{\mathbb{N}}
\def\R{\mathbb{R}}
\def\C{\mathbb{C}}

\hyphenation{imme-dia-te-ly}
\hyphenation{con-ti-nui-ty}
\hyphenation{se-ve-ral}
\hyphenation{e-qui-va-lent}
\hyphenation{e-qui-va-len-ces}
\hyphenation{fo-llo-wing}
\hyphenation{a-na-lo-gous-ly}

%%%%%%%%%%%

\title{The eigenvalues of limits of radial Toeplitz operators}
\author{Daniel Su\'{a}rez}

\date{\today}
\maketitle

%%  BEGIN ABSTRACT
\begin{quotation}
\noindent
\mbox{ } \hfill      {\sc Abstract}    \hfill \mbox{ } \\
 \hfill \mbox{ }\\
{\small \noindent
Let $A^2$ be the Bergman space on the unit disk.
A bounded operator $S$ on $A^2$ is called radial if $Sz^n = \lambda_n z^n$ for all $n\ge 0$,
where $\lambda_n$ is a bounded sequence of complex numbers.
We characterize the eigenvalues of radial operators that can be approximated
by Toeplitz operators with bounded symbols.
}
\end{quotation}
%% END ABSTRACT

\maketitle

\setcounter{section}{0}

\section{Introduction and preliminaries}

The Bergman space $\berg$ is the closed subspace of analytic functions in $L^2(\disc, dA)$, where
$\disc$ is the open unit disk and $dA$ is the normalized area measure.
The functions $e_n(z) := \sqrt{n+1}\, z^n$, with $n\geq 0$, form the standard orthonormal base of $\berg$.
We denote by $\L(\berg)$ the algebra of bounded operators on $\berg$. If $a\in \linf$,
the Toeplitz operator with symbol $a$ is
$$
T_af(z) := \int_D \frac{a(w)f(w)}{(1-\ov{w}z)^2}\, dA(w) , \ \ f\in\berg .
$$
It is immediate that $\|T_a\| \leq \|a\|_\infty$.
A Toeplitz operator $T_a$ is diagonal with respect to the standard base
(i.e.: $T_a e_n = \lambda_n e_n$ for some $\lambda_n \in \C$, $n\geq 0$) if and only if $a(z)=a(|z|)$.
By analogy, we say that $S\in\L(\berg)$ is a radial operator if it is diagonal with respect to the standard base.
Radial operators, mostly Toeplitz, have been studied by several authors
(see \cite{a-c}, \cite{g-v}, \cite{KoZh}  and \cite{zor}),
mainly because they are among the few Toeplitz operators on $\berg$ that we can reasonably understand so far.
Despite this fact, some central problems are still open.
Consider the space $\radop:=\{ T_b  : b\in\linf \mbox{ radial}\}$, and the Toeplitz algebra
$$
\toep := \mbox{ the closed subalgebra of $\L(\berg)$ generated by } \{ T_a : a\in \linf \} .
$$
In \cite{sua4} it is proved that any radial operator $S\in\toep$ can be approximated by
operators in $\radop$, whose symbols are constructed from $S$ in a canonical way (see Theorem \ref{oldmain} below).
Since radial operators are determined by their eigenvalues, two problems come readily to mind.
Can we characterize the eigenvalue sequences of operators in $\radop$?
What about the closure of $\radop$?
As we shall see immediately, the first question was settled more than 80 years ago. %almost 60 years ago.
The present paper deals with the second question.

If $b$ is a bounded radial function, the use of polar coordinates shows that its eigenvalue sequence
$\lambda(T_b)$ is given by
\begin{equation}\laba{eigg}
\lambda_n(T_b)  =  \la be_n,e_n \ra = (n+1)\int_0^1 b(r) \, r^{2n} \, 2r dr
= (n+1)\int_0^1 b(t^{1/2}) t^n \, dt .
\end{equation}
That is, a sequence $\{ \lambda_n \}$ forms the eigenvalues of a Toeplitz operator with bounded radial symbol
if and only if $\{ \lambda_n/(n+1) \}$ is the moment sequence of a bounded function on $[0,1]$.
In 1921 Hausdorff characterized the moment sequences of measures of bounded variation
on the interval $[0,1]$.

\begin{fed}
Let $m\geq 0$ be an integer and\/ $x=\{ x_n \}_{n\geq 0}$ be a sequence of complex numbers.
The $m$-difference of $x$, denoted $\Delta^m(x)$, is the sequence defined by
$$
\Delta^m_n x:= (-1)^m \sum_{j=0}^m   \left(^m_{\,j}\right)
(-1)^j \,x_{n+j} , \ \mbox{ for }\   n\geq 0,
$$
where ${m\choose{j}}= m!/(m-j)! j!$.
\end{fed}

\noi
Further elaboration of Hausdorff's moment theorem showed that
given a sequence $x$, there  exists a function $a\in L^\infty [0,1]$ such that
$\int_0^1 a(t) t^n \, dt = x_n$ for all $n\geq 0$
if and only if the expression $(k+1) \left(^{\,k}_m\right)  \left|\Delta^m_{k-m} \,x \right|$ is bounded
for all $0\leq m\leq k$ (see \cite[Ch.$\,$III]{wid}).
Together with \eqref{eigg}, this implies that a sequence $\lambda$ is formed by the eigenvalues of
some $T_b$, with $b\in L^\infty(\disc)$ radial,
if and only if there is a constant $C>0$ such that
\begin{equation}\laba{haus}
 (k+1)  \left(^{\,k}_m\right)  \left|\Delta^m_{k-m} \,\mu \right| \leq C
\,\ \mbox{ for all $\,0\leq m\leq k$, where }\,    \mu_n:= \frac{\lambda_n}{n+1} .
\end{equation}

\noi
Since $\|S\|=\|\lambda(S)\|_{\li}$ for any radial operator
with eigenvalue sequence $\lambda(S)$,
it is clear that $S$ is in the closure of\/ $\radop\,$
if and only if $\lambda(S)$
is in the $\li$-closure of the sequences that satisfy \eqref{haus}.
The obvious inconvenient with this characterization is that this property
is very hard to check.
Also, it is difficult to construct such sequences
without the a priori knowledge that the corresponding operator
is a limit of Toeplitz operators with bounded symbols.
We provide here two characterizations of these sequences that turn out to be much simpler than \eqref{haus}.
The resulting eigenvalues consist of the $\li$-closure of sequences $\lambda$ satisfying
any of the conditions:
$$
\sup_{n\geq 0} \, (n+1) \left|\Delta^{1}_n \,\lambda \right| < \infty
\ \ \mbox{ or }\ \
\sup_{n\geq 0} \, (n+1)^2 \left|\Delta^{2}_n \,\lambda \right| < \infty .
$$
This answers affirmatively a question that I posed in \cite{sua4}.
The precise relation between the above conditions and \eqref{haus} is made explicit in
Section \ref{appradto}.

\subsection{The $n$-Berezin transform}

\noi
If $n$ is a nonnegative integer and $z\in\disc$, consider the function
$$
K_z^n (\om) = \frac{1}{  (1-\ov{z}\om)^{2+n}  } \ \ (\om\in\disc) .
$$
When $n=0$ this function is the reproducing kernel for the space $\berg$.
\begin{fed}
The $n$-Berezin transform of an operator $S\in \L(\berg)$ is defined as
$$
B_n (S) (z) :=
(n+1) (1-|z|^2)^{2+n}  \sum_{j=0}^{n}  \left(^n_j\right)
(-1)^j  \,
\la S (\om^j K_z^n) , \om^j K_z^n \ra ,
$$
where $\la \, ,\ra$ is the usual integral pairing.
\end{fed}
\noi
It is not difficult to prove that $B_n(S)\in L^\infty (\disc)\cap C^\infty(\disc)$, and that there is a constant
$C(n)>0$ such that $\|B_n(S)\|_\infty \leq C(n) \|S\|$.
If $a\in L^\infty (\disc)$, the binomial expansion of $(1-|\om |^2 )^n$ and a conformal change of variables
yields
$$
B_n (a)(z)  := B_n (T_a)(z)= \int_D a( \varphi_z (\xi ) )
(n+1) (1-|\xi |^2 )^n \, dA(\xi ) ,
$$
where $\varphi_z(w) = (z-w)/(1-\ov{z}w)$ is the automorphism of the disk that interchanges $0$ and $z$.
That is, the above formula defines the $n$-Berezin transform of a function $a\in \linf$.
Since $(n+1) (1-|\xi |^2 )^n  dA(\xi )$ is a probability measure with total mass accumulating at $0$ when
$n\rr \infty$, it is clear that $B_n(a) \rr a$ pointwise if $a\in\linf$ is continuous. In particular,
\begin{equation}\laba{long1}
B_n (B_0(S)) \rr B_0(S)  \ \mbox{  for any $S\in \L(\berg)$.}
\end{equation}
The (conformally) invariant Laplacian is $\tilde{\Delta} = (1-|z|^2)^2 \Delta$, where
$\Delta=\partial \ov{\partial}$ is a quarter of the standard Laplacian, and
$\partial, \, \ov{\partial}$ are the traditional Cauchy-Riemann operators.
It is easy to check that
$(\deltab f)\circ \psi = \deltab (f\circ \psi)$ for every $f\in C^2(\disc)$ and $\psi\in \mbox{Aut}(\disc)$.

We summarize next some of the
properties of $B_n$ that will be used in the paper.
The proofs are in \cite{sua3}.
Let $S\in \L (\berg)$ and\/ $n\geq 0$. Then

\begin{equation}   \laba{deltu}
\tilde{\Delta} B_n (S) = (n+1)(n+2) (B_n(S) - B_{n+1}(S)) .
\end{equation} \vspace{-2mm}
\begin{equation}   \laba{commu}
(B_k B_j) (S) = (B_j B_k) (S) \mbox{ for all $j,k \geq 0$} .
\end{equation} %\vspace{-2mm}
Observe that \eqref{deltu} implies that $\deltab B_n(S) \in \linf$ for any $S\in\L(\berg)$, which allows us
to further apply $B_k$ to this function for any $k\geq 0$.
It follows immediately from \eqref{deltu} and \eqref{commu} that
\begin{equation}   \laba{deltoide}
\deltab B_k (B_nS)= B_k \deltab (B_nS).
\end{equation}
Finally, it is easy to prove that if $S\in\L(\berg)$ is radial, so is the function $B_n(S)$.

\section{Two sequence spaces}\label{twoseq}
\equnew

\noi
Let $\ell^\infty$ be the Banach space of bounded complex sequences indexed from $n\geq 0$.

\begin{fed}
Consider the linear subspaces of\/ $\li$:
\begin{align*}
d_1 &:=\left\{x\in\li:\  \|x\|_{d_1} = \sup_{n}\, (n+1) |\Delta^1_n(x)|<\infty\right\}   \ \mbox{and} \\
d_2 &:=\left\{x\in\li:\  \|x\|_{d_2} = \sup_{n}\, (n+2)^2 |\Delta^2_{n}(x)|<\infty\right\}.
\end{align*}
\end{fed}

\noi Observe that $\|\ \|_{d_1}$ and $\|\ \|_{d_2}$ are semi-norms that vanish only at constant sequences.
\begin{lemma}\laba{estic3}
Given $C>4$, for every $n\geq C$ there exists $r\in [1,4]$ such that\/ $m:=\frac{C}{C-r}n$ %\, n\in\N$
is the unique integer that satisfies
\begin{equation}\laba{dos}
\sum_{k=n+1}^{m}\frac{C}{k^2}\leq \frac{1}{n}\peso{and} \sum_{k=n+1}^{m+1}\frac{C}{k^2}> \frac{1}{n}.
\end{equation}
\end{lemma}
\bdem
Fix $C> 4$ and suppose that $n\geq C$. For $m\ge n+1$, we have
$$
\frac{1}{2} \left( \frac{1}{n}-\frac{1}{m} \right) \le
\frac{1}{n+1}-\frac{1}{m+1} = \int_{n+1}^{m+1}\frac{1}{x^2}\ dx
\leq
\sum_{k=n+1}^{m}\frac{1}{k^2}
\leq
\int_{n}^{m}\frac{1}{x^2}\ dx = \frac{1}{n}-\frac{1}{m}.
$$
Straightforward estimates from these inequalities show that

$$
C\sum_{k=n+1}^{m}\frac{1}{k^2}\leq \frac{1}{n} \peso{if}  n+1\leq m\leq \frac{C}{C-1}\,n,
$$
and %%straightforward estimates from the first inequality show that
$$
C\sum_{k=n+1}^{m}\frac{1}{k^2}\geq \frac{2}{n} \peso{if}  m\geq \frac{C}{C-4}n.
$$

\noi Hence, there is $m\in\N$ between $\frac{C}{C-1}\,n$ and $\frac{C}{C-4}\,n$ satisfying \eqref{dos}.
Since the function $f(x)=\frac{C}{C-x}n$ is continuous on $[1,4]$, the mean value theorem gives
$r\in [1,4]$ such that $m= \frac{C}{C-r}n$.
\edem

\begin{lemma}\laba{estic4}
Given $\eps>0$, there exists $C=C(\eps)>4$ large enough so that for all $n\geq C$,
\begin{equation}\laba{ab1}
E(C,n):=\left[\frac{1}{n}-\frac{C}{(n+1)^2}\right]+\ldots+\left[\frac{1}{n}-\sum_{k=n+1}^m\frac{C}{k^2}\right]
< \eps
\end{equation}
and
\begin{equation}\laba{ab2}
\sum_{k=n+1}^{m}\frac{1}{k} <\eps,
\end{equation}
where $m=m(C,n)$ is the integer  satisfying \eqref{dos}.
\end{lemma}
\bdem
First observe that
$$
E(C,n) =\frac{m-n}{n}-C\sum_{k=n+1}^m \frac{(m-k+1)}{k^2}\leq \frac{m-n}{n}.
$$
Also, since $\log(1+x) \leq x$ for $x\geq 0$,
$$
\sum_{k=n+1}^{m}\frac{1}{k}\leq \int_{n}^{m}\frac{1}{x}\ dx %%=\left.\log(x)\right|_{n}^m
=\log\left(\frac{m}{n}\right) \leq \frac{m}{n}-1. %\left(\frac{m}{n}-1\right).
$$
By Lemma \ref{estic3} there is $r=r(C,n)\in [1,4]$ such that $m=\frac{C}{C-r}\, n$.
Hence, replacing $m$ by this expression we get
$$
\frac{m}{n}-1 = \frac{C}{C-r}-1 \leq \frac{C}{C-4}-1 = \frac{4}{C-4} <\eps
$$
if $C$ is large enough.
\edem

\noi I learnt the proof of the next proposition from Jorge Antezana (personal communication), to whom I
am grateful.

\begin{pro}\laba{drury}
The following statements hold
\begin{enumerate}
\item[{\em (1)}] If\/ $x\in d_2$ then $\| x\|_{d_1} \leq \| x \|_{d_2}$.
\item[{\em (2)}] $\clau{d_1}=\clau{d_2}$.
\end{enumerate}
\end{pro}
\bdem
Let $x \in d_2$. Fix $j\geq 0$ and let $n\geq j$. Then
\begin{eqnarray}\laba{dedos}
|\Delta^1_{n+1}(x)-\Delta^1_j(x)|
&\leq&   \sum_{k=j}^n |\Delta^1_{k+1}(x) -  \Delta^1_{k}(x)|  \ = \  \sum_{k=j}^n |\Delta^2_{k}(x)| \nonumber \\
&\leq&   \sum_{k=j}^n \frac{\| x \|_{d_2}}{(k+1)(k+2)}
\ = \ \| x \|_{d_2} \left( \frac{1}{j+1} -  \frac{1}{n+2}\right) \! .
\end{eqnarray}
Hence, $\Delta^1_{n+1}(x)$ is a Cauchy sequence, and since the sequence $x$ is bounded,
$\Delta^1_{n+1}(x) \rr 0$.
Taking limit when $n\rr \infty$ in \eqref{dedos} we obtain
$|\Delta^1_j(x)| \leq  \frac{\| x \|_{d_2}}{j+1}$ for $j\geq 0$.
This proves (1).
In particular, $d_2\subseteq d_1$, and the proof of (2) is reduced to see that
$d_1\subset \clau{d_2}$.
So, let $x\in d_1$ and $\eps>0$. We can assume without loss of generality that
$\|x\|_{d_1}\leq 1$ and that $x_n\in\R$ for every $n\geq 0$.
Pick $C=C(\eps)>4$ as in Lemma \ref{estic4}
and define $y\in d_2$ as:
$$
y_n = \left\{
\begin{array}{ll}
x_n                                            & \mbox{if $\, n\leq \max\{ 2/\eps, \, C\}$} \\*[1mm]
y_{n-1}+\delta                                 & \mbox{if $\, n>\max\{ 2/\eps, \, C\}$,}\\
\end{array}
\right.
$$
$$
\mbox{where $\delta$ minimizes}\ \,  |(y_{n-1}+\delta)-x_{n}|
\mbox{ under the restrictions }
 |\delta|  \leq \frac{1}{n}, \ |(y_{n-1}-y_{n-2})-\delta| \leq \frac{C}{n^2}  .\\
$$

\noi
We aim to prove that $\|x-y\|_\infty \leq 5\eps$.
First observe that given $n_0$ such that $|x_{n_0}-y_{n_0}|<\eps$,
it is enough to estimate $|x_n-y_n|$ for
all the subsequent values of $n$ until the first time that $\sgn (x_n-y_n) \neq \sgn (x_{n+1}-y_{n+1})$,
because this change of sign implies that $|x_{n+1}-y_{n+1}| \leq 2/(n+1)  <\eps$.
So, suppose that $n_0$ is such that $x_{n_0} < y_{n_0}$ (the analysis is symmetrical for $x_{n_0} > y_{n_0}$),
and let $n_1$ be the first integer $>n_0$ such that $x_{n_1} \geq y_{n_1}$.
We estimate how much $y_n-x_n$ can grow for $n_0< n< n_1$.
Since $x_n<y_n$ for those values of $n$,
\bequ\laba{urxe}
y_{n}-y_{n-1}= y_{n-1}-y_{n-2} -\frac{C}{n^2} =(y_{n_0}-y_{n_0-1}) -\! \!\sum_{k=n_0+1}^n \frac{C}{k^2}
%%\ \mbox{ for all }\   n_0 <n< n_1.
\eequ
for all $n_0 <n< n_1$.
Consider two cases.

Case 1: $y_{n_0}-y_{n_0-1}\leq 0$. It follows from \eqref{urxe} that for $n_0 <n< n_1$,
\begin{eqnarray*}%\laba{goz}
y_n-x_n &=&  (y_{n_0} -x_{n_0})  + (y_n-y_{n_0}) - (x_n -x_{n_0})  \\
&\le& (y_{n_0} -x_{n_0})
 - \left[ \frac{C}{(n_0+1)^2} + \cdots +  \sum_{k=n_0+1}^n \frac{C}{k^2}\right]
 -\sum_{j=n_0+1}^n (x_{j} -x_{j-1})   \\
&\le&
(y_{n_0} -x_{n_0})
+ \left[ |x_{n_0+1} -x_{n_0}| - \frac{C}{(n_0+1)^2} \right] + \cdots +
\left[ |x_{n} -x_{n-1}| - \sum_{k=n_0+1}^n \frac{C}{k^2} \right]
\end{eqnarray*}
Let $m_0>n_0$ be the integer associated with $n_0$ by \eqref{dos}.
If $n>m_0$, \eqref{dos} and $|x_{n} -x_{n-1}| \leq \frac{1}{n} \leq \frac{1}{n_0}$ imply that
the corresponding summand in square brackets must be negative.
Thus,
\begin{align*}
y_n-x_n
&\leq
(y_{n_0} -x_{n_0})
+ \left[ \frac{1}{n_0}  - \frac{C}{(n_0+1)^2} \right] + \cdots +
\left[ \frac{1}{n_0} - \sum_{k=n_0+1}^{m_0} \frac{C}{k^2} \right]\\
&= (y_{n_0} -x_{n_0}) + E(C,n_0) \leq (y_{n_0} -x_{n_0}) + \varepsilon ,
\end{align*}
where the las inequality comes from \eqref{ab1}.

Case 2: $y_{n_0}-y_{n_0-1}> 0$.
If $n_0 \leq n< n_1$ is any integer such that
\begin{equation}\laba{mixi}
y_{k}-y_{k-1}> 0  \ \mbox{ for }\   k= n_0, \ldots , n,
\end{equation}
then \eqref{urxe} implies that
$$
0< y_n-y_{n-1} =  (y_{n_0}-y_{n_0-1}) - \sum_{k=n_0+1}^n \frac{C}{k^2}
\le \frac{1}{n_0} - \sum_{k=n_0+1}^n \frac{C}{k^2},
$$
which together with the definition of $m_0$ forces $n \leq m_0$.
So, \eqref{ab1} and \eqref{ab2} give  %%going to \eqref{goz}, we see that
\begin{align}\label{qlpa}
y_n-x_n &=  (y_{n_0} -x_{n_0})  + (y_n-y_{n_0}) - (x_n -x_{n_0}) \nonumber \\
&\leq  (y_{n_0} -x_{n_0})
+ \left[ \frac{1}{n_0} - \frac{C}{(n_0+1)^2} \right] + \cdots +
\left[  \frac{1}{n_0} - \sum_{k=n_0+1}^{m_0} \frac{C}{k^2} \right]
+ \sum_{k=n_0+1}^{m_0} |x_{k} -x_{k-1}| \nonumber \\
&\leq  (y_{n_0} -x_{n_0}) + E(C,n_0) + \sum_{k=n_0+1}^{m_0} \frac{1}{k}
\leq  (y_{n_0} -x_{n_0}) + 2\varepsilon .
\end{align}
If $n$ is the largest integer satisfying \eqref{mixi}, then
either $n+1= n_1$ (and we are done) or $n+1$ is in Case 1, meaning that $y_{n+1}-y_n \leq 0$
(while $y_{n+1}>x_{n+1}$).
Hence, the estimate of Case 1 and \eqref{qlpa} show that for all $n+1 \leq k <n_1$,
$$
y_k-x_k \leq (y_{n+1}-x_{n+1}) + \eps
\leq \frac{2}{(n+1)} + (y_{n}-x_{n}) + \eps
< \eps + (y_{n_0} -x_{n_0}) + 2\varepsilon + \eps .
$$
That is, we have shown that $y_k-x_k \leq y_{n_0} -x_{n_0} +4\eps$ for all $n_0\leq k< n_1$.
By the symmetry of the case $x_{n_0} > y_{n_0}$ and the comments that follow
the definition of $y$, we get $\|y-x\|_\infty < 5\eps$.
\edem

\section{The invariant Laplacian of an operator}
\equnew

\begin{fed}
Let
$$\cald =\{ S\in \L(\berg):  \exists \,
T\in \L(\berg) \,\mbox{ such that }\, \deltab B_0(S) = B_0(T) \} ,
$$
and define $\deltab: \cald \, \rr \, \L(\berg)$ by\/ $\deltab S = T$.
\end{fed}

\begin{lemma}\laba{pobo}
If $S_n, S\in\L(\berg)$, with $S_n \rr S$ in the weak operator topology. Then
$$
B_0 (S_n) \rr  B_0 (S)   \ \mbox{ and }\  \deltab B_0 (S_n) \rr \deltab B_0 (S) \ \mbox{ pointwise.}
$$
\end{lemma}
\bdem
We only prove the assertion for $\deltab B_0$, since the proofs are analogous.
It is clear that if $k$ is a non negative integer,
$\ov{\partial}_z (\ov{z}^k K_z)(w)$ is a bounded analytic function
of $w$. Thus,
$$
\Delta |z|^{2k} \la S_n K_z,K_z\ra =  \la S_n
\,\ov{\partial}_z(\ov{z}^k K_z), \ov{\partial}_z(\ov{z}^k K_z)\ra
\,\rr\, \la S\,\ov{\partial}_z(\ov{z}^k K_z),
\ov{\partial}_z(\ov{z}^k K_z)\ra = \Delta |z|^{2k} \la SK_z,K_z\ra ,
$$
with point convergence on $z$. In particular,
$$
\deltab B_0 (S_n) = (1-|z|^2)^2 \Delta  [(1+|z|^4-2|z|^2)\la S_n
K_z,K_z\ra]
$$
converges pointwise to $\deltab B_0 (S)$.
\edem

\begin{lemma}\laba{lamga}
For $\lambda\in d_2$ consider the sequence $\gamma$ given by\/
\begin{equation}\laba{gamy}
\gamma_{n} :=
\left\{
\begin{array}{ll}
2(\lambda_1-\lambda_0),                                                                 & \mbox{if $\, n=0$} \\*[1mm]
(n+1)\, [(n+2) (\lambda_{n+1} -  \lambda_{n})  -  n (\lambda_{n} -  \lambda_{n-1}) ],   & \mbox{if $\, n\geq 1$}
\end{array}
\right.
\end{equation}
Then
$$
6^{-1}\,  \|\lambda\|_{d_2} \leq       \| \gamma\|_\infty        \leq 6\, \|\lambda\|_{d_2}.
$$
\end{lemma}
\bdem
Setting $\lambda_{-1}=0$, for $n\geq 0$ we have
\begin{eqnarray*}
\gamma_{n}
&=& (n+2)(n+1) \, \Delta^1_{n}(\lambda)  -  (n+1)n \, \Delta^1_{n-1}(\lambda) \\
&=& (n+2) \, b_{n+1}  -  (n+1)\, b_{n} ,         %\ \ \mbox{ for $n\geq 0$,}
\end{eqnarray*}
where $b_n:=n  \Delta^1_{n-1}(\lambda)$.
Therefore
$$
(n+2)\, |b_{n+1}| =
\left|(n+2) b_{n+1} -   1 b_{0} \right|=
\left| \sum_{j=0}^n [(j+2) \, b_{j+1}  -  (j+1)\, b_{j} ] \right|
\leq (n+1) \| \gamma\|_\infty ,
$$
leading to $(n+1)|\Delta^1_{n} (\lambda) |=|b_{n+1}| \leq  \| \gamma\|_\infty$, for $n \geq 0$.
That is, $\|\lambda\|_{d_1} \leq \| \gamma\|_\infty$.
On the other hand, if $n\geq 1$,
\begin{eqnarray*}
\gamma_{n}
&=&    (n+1)\, [(n+2) (\lambda_{n+1} -  2\lambda_{n} + \lambda_{n-1}) +  2(\lambda_{n} -  \lambda_{n-1}) ] \\*[0.7mm]
&=&    \left(\frac{n+2}{n+1}\right) (n+1)^2 \Delta^2_{n-1}(\lambda)
+  2\left(\frac{n+1}{n}\right)n\Delta^1_{n-1}(\lambda) .
\end{eqnarray*}
Hence,
$$
\| \gamma \|_{\li}
\leq   2 \| \lambda \|_{d_2} + 4 \| \lambda \|_{d_1}
\stackrel{\mbox{\scriptsize{by Prop.$\,$\ref{drury}}}}{ \leq }     6 \| \lambda \|_{d_2},
$$
and since
$$
|(n+1)^2 \Delta^2_{n-1}(\lambda)| =
\left|\frac{(n+1)}{(n+2)} \, \gamma_{n}    -   2\frac{(n+1)^2}{(n+2)n}\,  n\Delta^1_{n-1}(\lambda)\right|
\leq
|\gamma_{n}| + 4n|\Delta^1_{n-1}(\lambda)|,
$$
then
$$
\| \lambda \|_{d_2} \leq         \| \gamma \|_{\li} +  4\| \lambda \|_{d_1}
\leq 5 \| \gamma \|_{\li} .
$$
\edem

\noi
The orthogonal projection onto the subspace generated by $e_n$ is $E_n f = \la f, e_n\ra e_n$,
where $n\geq 0$ and $f\in\berg$.
Thus, a bounded operator $S$ is radial if and only if it can be written as %% admits the decomposition
$S=\sum_{n\geq 0} \lambda_n E_n$, where $\lambda\in\li$ is the sequence of its eigenvalues.
Also, observe that the reproducing property of $K_z^0$ shows that

\begin{equation}\laba{berto}
(1-|z|^2)^2 B_0(E_n)(z) = (1-|z|^2)^2  \la K_z^0, e_n \ra    \la e_n, K_z^0 \ra =  (1-|z|^2)^2 |e_n(z)|^2.
\end{equation}

\begin{pro}\laba{radlap}
Let $S\in\L(\berg)$ be a radial operator with eigenvalue sequence $\lambda$.
Then $S\in \cald$ if and only if\/ $\lambda\in d_2$, in which case
\begin{equation}\laba{dlin}
\deltab \sum_{n\geq 0} \lambda_n E_n = \sum_{n\geq 0} \gamma_{n}  E_n,
\end{equation}
where $\gamma$ is given by \eqref{gamy}.
Thus, $6^{-1}\,  \|\lambda\|_{d_2} \leq       \| \deltab S\|  \leq 6\, \|\lambda\|_{d_2}$.
\end{pro}
\bdem
Since the partial sums of $\sum \lambda_n E_n$ tend to $S$ in the strong operator topology,
Lemma \ref{pobo} implies that
$$
\deltab B_0 \left( \sum \lambda_n E_n \right) =\sum \lambda_n \deltab B_0 (E_n).
$$
By \eqref{berto},
\begin{align*}
\deltab B_0 (E_n)(z) &= (n+1)(1-|z|^2)^2 \Delta (1-|z|^2)^2 |z|^{2n} \\*[1mm]
&= (n+1)(1-|z|^2)^2 (n^2 |z^{n-1}|^2 + (n+2)^2 |z^{n+1}|^2-2(n+1)^2|z^{n}|^2) .
\end{align*}
Then
\begin{align*}
\deltab B_0 (S)
&= (1-|z|^2)^2  \sum_{n\geq 0} \lambda_n  (n+1) (n^2 |z|^{2(n-1)} + (n+2)^2 |z|^{2(n+1)}-2(n+1)^2|z|^{2n}) \\
&= (1-|z|^2)^2  \sum_{n\geq 0} |z|^{2n} (n+1)^2 [(n+2)\lambda_{n+1} + n\lambda_{n-1} -2(n+1)\lambda_n ] \\
&= (1-|z|^2)^2  \sum_{n\geq 0} |e_n(z)|^{2} (n+1) [(n+2)\lambda_{n+1} + n\lambda_{n-1} -2(n+1)\lambda_n ] ,
\end{align*}
where we are taking $\lambda_{-1}=0$, and the second equality comes from regrouping the series, which is absolutely
and uniformly convergent on compact sets of $\disc$.
That is,
$$\deltab B_0 (S) = (1-|z|^2)^2  \sum_{n\geq 0} \gamma_n |e_n(z)|^{2},$$
where
$\gamma_{n} =    (n+1)\, [  (n+2)\lambda_{n+1}+n  \lambda_{n-1}-2 (n+1)  \lambda_{n}  ]$.
If $\lambda\in d_2$, Lemma \ref{lamga} says that
$\gamma\in \li$. So, the operator $T:= \sum_{n\geq 0} \gamma_n E_n$ is bounded, and \eqref{berto} with
Lemma \ref{pobo} imply that
$$
B_0(T)= (1-|z|^2)^2  \sum_{n\geq 0} \gamma_n |e_n(z)|^{2}= \deltab B_0(S).
$$
Reciprocally, suppose that $T$ is a bounded operator that satisfies $B_0(T)= \deltab B_0(S)$.
Writing $K^{0}_z(w) = \sum\ov{e_m(z)} e_m(w)$ we get
$$
B_0(T)(z) = (1-|z|^2)^2\la TK^{0}_z,K^{0}_z\ra
= (1-|z|^2)^2  \sum_{n,\,m=0}^\infty  \la T e_n, e_m\ra \,  \ov{e_n(z)} e_m(z),
$$
which clearly implies that $\la T e_n, e_m\ra =0$ for $n\neq m$ and $\la T e_n, e_n\ra =\gamma_n$.
Therefore, $\gamma\in\li$ and Lemma \ref{lamga}  implies that $\lambda\in d_2$.

In either case, $\deltab S = \sum \gamma_n (e_n\otimes e_n)$, which proves \eqref{dlin}, and since
$\|\deltab S\| = \|\gamma\|_{\li}$, the last assertion of the proposition follows from Lemma \ref{lamga}.
\edem

\section{Approximation by radial Toeplitz operators}\label{appradto}
\equnew

\begin{lemma}\laba{crirad}
Suppose that $S\in \L(\berg)$ is such that
$\|T_{\deltab B_k(S)}\| \leq\/ C$ independently of\/ $k$. Then $T_{B_k(S)} \rr S$.
\end{lemma}
\bdem  By \eqref{deltu},
$
T_{\deltab B_k (S)} = (k+1)(k+2) (T_{B_k (S)} - T_{B_{k+1} (S)}) .
$
So,
$$
T_{B_0 (S)} -   \sum_{k=0}^m \frac{T_{ \deltab B_k(S) }}{(k+1)(k+2)}
= T_{B_{m+1} (S)} ,
$$
and since $\|T_{ \deltab B_k(S) }\| \leq C$, the series of the norms is convergent, which implies
the convergence of $T_{B_{m} (S)}$, say to $R\in \L(\berg)$.
Since $B_0$ is a bounded operator from $\L(\berg)$ in $L^\infty$, we also have that
$B_0(T_{B_{m} (S)}) \rr B_0(R)$ in $L^\infty$-norm.
On the other hand, \eqref{commu} and \eqref{long1} imply
$$
B_0( T_{B_{m} (S)} ) = B_0 B_m (S) = B_m B_0 (S) \rr B_0(S) \ \mbox{ pointwise.}
$$
This means that $B_0(S)=B_0(R)$, and since $B_0$ is one-to-one, $S=R$.
\edem

\noi
We recall that the Toeplitz algebra, $\toep$, is formed by all the operators that can be approximated by polynomials
of Toeplitz operators with bounded symbols. The two results in the following theorem are
Corollary 3.2 and Theorem 3.3 of \cite{sua4}, respectively.

\begin{theo}\laba{oldmain}
Let $S\in\L(\berg)$ be a radial operator. Then
\begin{enumerate}
\item[{\em 1.}] $\|T_{B_k(S)}\| \leq \|S\|$.
\item[{\em 2.}] $S\in\toep$ if and only if\/ $T_{B_k(S)} \rr S$.
\end{enumerate}
\end{theo}

\noi It is easy now to finish the proof of the main result in this paper.

\begin{theo}\laba{newmain}
Let $S\in \L(\berg)$ be a radial operator with eigenvalue sequence $\lambda(S)$. Then the following
statements are equivalent
\begin{enumerate}
\item[{\em (1)}]  $S\in \toep\,$ {\em (}or equivalently, $T_{B_k(S)}\rr S${\em )}\vspace{-2mm}
\item[{\em (2)}]  $\lambda(S)   \in \ov{d}^{\,\ell^\infty}_2$ \vspace{-2mm}
\item[{\em (3)}]  $\lambda(S)  \in \ov{d}^{\,\ell^\infty}_1$
\end{enumerate}
\end{theo}
\bdem
Observe that Proposition \ref{drury} gives the equivalence between (2) and (3).
We shall prove that (1) is equivalent to (2).
It is quite easy to show that if $b$ is a bounded radial function then its eigenvalue sequence
$\lambda(T_b)$ is in $d_2$. Indeed, if $n\geq 1$, \eqref{eigg} yields
$$
|\Delta_{n-1}^2 (\lambda(T_b))|
\leq  \int_0^1 |b(t^{1/2})|\,  |(n+2)t^{n+1}-2(n+1)t^{n}+nt^{n-1}| \, dt
\leq     \frac{8\,\|b\|_\infty}{(n+2)^2}  .
$$
If (1) holds then
$$\lambda(T_{B_k(S)}) \stackrel{\ell^\infty}{\rr} \lambda(S)
\peso{when} n\rr\infty.$$
So, $\lambda(S) \in \ov{d}^{\,\ell^\infty}_2\!$.
Now suppose that $\lambda(S)$ is the $\ell^\infty$-limit of a sequence $\lambda_j$ contained in $d_2$
(here $\lambda_j$ denotes the whole sequence, not the $j$-entry of a sequence). If $S_j$ is the radial operator
with eigenvalues $\lambda(S_j)= \lambda_j$, then $S_j\rr S$ in $\L(\berg)$-norm.
If we show that $T_{B_k(S_j)}\rr S_j$ when $k\rr\infty$ for every fixed value of $j$,
then $S\in\toep$ and (1) will follow. That is, we can assume that $\lambda(S)\in d_2$.
By Proposition \ref{radlap} then $S\in \cald$ and $\|\deltab S\| \leq 6 \|\lambda(S)\|_{d_2}$.
Since $\deltab S$ is a radial operator, Theorem \ref{oldmain} says that
$\|T_{B_k(\deltab S)}\| \leq \|\deltab S\|$.
Furthermore, by \eqref{commu} and \eqref{deltoide},
$$
B_0 \deltab B_k( S)=  \deltab B_0 B_k(S) = \deltab B_k B_0 (S) = B_k \deltab  B_0 (S)
= B_k   B_0 (\deltab S) =  B_0 B_k  (\deltab S),
$$
and since $B_0$ is one-to-one, $\deltab B_k( S)= B_k  (\deltab S)$.
Putting all this together gives
$$
\|T_{\deltab B_k(S)}\| =  \|T_{B_k(\deltab S)}\| \leq \|\deltab S\| \leq 6 \|\lambda(S)\|_{d_2}
$$
for all $k$. Lemma \ref{crirad} then says that  $T_{B_k(S)}\rr S$.
\edem

\noi
A direct comparison between the conditions defining $d_1$ and $d_2$ with \eqref{haus} shows that a sequence
$\lambda$ satisfies \eqref{haus} for
\begin{align*}
m& =0 \piso{and} k\geq 0              \ \  \Leftrightarrow  \ \    \lambda\in\li ,  \\
m& =0, 1 \piso{and} k\geq m        \ \  \Leftrightarrow  \ \   \lambda\in d_1 ,\\
m& =0, 1, 2 \piso{and} k\geq m   \ \  \Leftrightarrow  \ \   \lambda\in d_2 .
\end{align*}
Therefore, if for any integer $p\geq 1$ we  define
$$
d_p :=\left\{x\in\C^{\N_0}:\,  x \ \, \mbox{satisfies \eqref{haus}$\ $ for $m=0, \ldots , p\ $ and $\ k\geq m$}
\right\} ,
$$
then $d_{p+1}\subset d_p \subset \li$, and the comment that follows \eqref{haus} together
with Theorem \ref{newmain} yield
$$
\ov{ \bigcap_{p\geq 1} d_p}^{\,\ell^\infty} \! \! = \, \ov{d}^{\,\ell^\infty}_1   \!   \! \! .
$$
In particular, an immediate consequence is the second assertion of Proposition \ref{drury}.
However, the assertion should be proved independently of this equality in order to avoid a cyclic argument.

Next we see two applications of the theorem.
Formula \eqref{eigg} defines a sequence $\lambda(b)$ for any radial function $b\in L^1(\disc)$, with
$$
\lambda_n(b)  =  (n+1)\int_0^1 b(t^{1/2}) t^n  dt,  \ \mbox{ for $n\ge 0$.}
$$
So, $b$ induces a bounded Toeplitz operator $T_b$ on $\berg$ if and only if
the sequence $\lambda(b)$ is bounded, with $\|T_b\|=\| \lambda(b)\|_{\li}$.
To this writing I do not know any geometric necessary and sufficient condition on $b$ for this to hold.
However, there is a well-known sufficient condition:
\begin{equation}\laba{lcon}
\left| \int_t^1 b(x^{1/2}) dx \right| \le C (1-t)  \ \mbox{ for all\/ $t\in[0,1]$},
\end{equation}
which turns out to be necessary when $b\ge 0$. Actually, when $b\ge 0$, \eqref{lcon} is a particular case
of a more general situation involving Carleson measures for Bergman spaces.
The next corollary shows that if $b$ satisfies \eqref{lcon} then $T_b$ is not only bounded, but it belongs to the
Toeplitz algebra $\toep$, and even to $\cald$.

\begin{coro}
Let $b\in L^1(\disc)$ be a radial function satisfying \eqref{lcon}.
Then $\|\lambda(b) \|_{\li} \le C$ and\/  $\|\lambda(b) \|_{d_2} \le 10\,C$.
In particular, $T_b\in\cald$ and hence, in $\toep$ {\em (}by Prop.$\ $\ref{radlap} and Thm.$\ $\ref{newmain}{\em )}.
\end{coro}
\bdem
For $n\ge 1$, integration by parts gives
$$
\lambda_n(b)
=   \int_0^1 \left[ \int_t^1 b(x^{1/2}) dx \right] (n+1)nt^{n-1} \, dt .
$$
Using \eqref{lcon} we immeditely see that
$|\lambda_n(b)| \le C$ for $n\ge 1$. For $n\ge 2\,$:
\begin{align*}
|\Delta^2_{n-1}(\lambda(b))|
&\le C \int_0^1 (1-t)t^{n-2} |(n+2)(n+1) t^2 - 2(n+1)n t  +  n (n-1)|  \, dt \\
&= C \int_0^1 (1-t)t^{n-2} |(n+1)n (1-t)^2 +2n (t^2-1) +2t^2|  \, dt \\
&\le 2C \int_0^1 (1-t)t^{n-2} [n^2 (1-t)^2  + t^2]  \, dt \\*[1mm]
&= 2C \left[    n^2 \, \frac{3! (n-2)! }{(n+2)!}
+ \, \frac{ n! }{(n+2)!} \right] \ \le \  \frac{10\, C}{(n+1)^2},
\end{align*}
where the last equality comes from $\int_0^1 (1-t)^p \, t^q \,dt =  p! \, q!/(p+q+1)!$
for integers $p, q\ge 0$.
Since $|\lambda_{0}(b)| = |\int_0^1 b(x^{1/2})\, dx|\le C$ by \eqref{lcon},
and $|\Delta^2_{0}(\lambda)|\le  |\lambda_2|+2|\lambda_1|+|\lambda_0|\le 3C$, the corollary follows.
\edem

\noi
It is known that if $S\in\L(\berg)$ is diagonal, then its essential spectrum $\sigma_e(S)$ is
formed by the limit points of its eigenvalues. In particular, since $\Delta^1_n(\lambda(b))\rr 0$
for any radial $b\in L^1(\disc)$, then $\sigma_e(T_b)$ is connected whenever $T_b$ is bounded.
Since also $\Delta^1_n(\lambda)\rr 0$ when $\lambda$ belongs to the $\li$-closure of $d_1$,
Theorem \ref{newmain} implies that $\sigma_e(S)$ is connected for every radial $S\in\toep$.
In \cite[Coro.$\,$2.10]{g-v}, Grudsky and Vasilevski show
examples of compact sets that can be the essential spectrum of $T_b$, for $b\in\linf$ radial.
We finish this paper by showing that if instead of $\radop$ %this class of operators
we take its closure, any nonempty, compact, connected set
is the essential spectrum of some operator in this class.

\begin{coro}
Let $E\subset \C$ be a nonempty, compact, connected set. Then there is a radial operator $S\in \toep$
such that $\sigma_e(S)=E$.
\end{coro}
\bdem
It is easy to construct a sequence $\lambda\in d_1$ whose limit points are exactly
the points of $E$. If $S$ is the radial operator with eigenvalue sequence $\lambda$, then $\sigma_e(S) =E$,
and Theorem \ref{newmain} says that $S\in\toep$.
\edem

\medskip

\noi {\bf Acknowledgements:}
The author is partially supported by the Ram\'on y Cajal program and the grants MTM2005-00544 and 2005SGR00774,
from the State Secretary of Education and Universities, Spain.

\newcommand{\foo}{\footnotesize}
\bigskip

 \noindent Daniel Su\'{a}rez\\
 Departament de Matem\`{a}tiques \\
 Universitat Aut\`{o}noma de Barcelona \\
 08193, Bellaterra, Barcelona \\
 Spain\\
\vspace{0.5mm} \noindent $\! \!${\foo dsuarez@mat.uab.es}

\end{document}